\documentclass[11pt]{amsart}

\usepackage{amsmath}
\usepackage{amssymb}
\usepackage{graphicx}
\usepackage{tikz-cd}

\newcommand{\St}{\operatorname{St}}
\newcommand{\mcG}{\mathcal{G}}

\newcommand{\mcR}{\mathcal{R}}
\newcommand{\Exp}{\operatorname{ EXP }}

\newcommand{\mcF}{\mathcal{F}}
\newcommand{\mcT}{\mathcal{T}}
\newcommand{\mcC}{\mathcal{C}}

\newcommand{\mcO}{\mathcal{O}}

\newcommand{\mcS}{\mathcal{S}}
\newcommand{\mcH}{\mathcal{H}}

\newcommand{\mcV}{\mathcal{V}}

\newcommand{\mcK}{\mathcal{K}}

\newcommand{\mcQ}{\mathcal{Q}}

\newcommand{\one}{\mathbf{1}}
\newcommand{\RR}{\mathbf{R}}

\newcommand{\q}{{q}}

\newcommand{\om}{\omega}
\newcommand{\Om}{\Omega}

\renewcommand{\sp}{\operatorname{sp}_-}
\newcommand{\CC}{\mcC^\infty(\Omega)}

\DeclareMathOperator{\ord}{ord}

\newcommand{\PG}{PGL(2,\mcQ_q)}

\newcommand{\dge}{d(o,g\cdot o)}
\newcommand{\Kk}{\mcK(\Omega)}

\newcommand{\Aut}{\mathcal{A}ut(\mcT)}

\newcommand{\inv}[1]{{#1}^{-1}}
\newcommand{\Norm}[1]{{\Vert #1 \Vert}}

\newcommand{\norm}[1]{\|{#1}\|}

\newcommand{\st}{\,;\;} 
\newcommand{\period}{\; .}    

\title[$q$-adic canonical states]{Representations of currents\\
taking values in $PGL(2,\Bbb Q_q)$}

\author[M.G. Kuhn]{Maria Gabriella Kuhn }
\address
{Dipartimento di Matematica \\
Universit\`a di Milano ``Bicocca'' \\
Via Cozzi 53 Ed.U5 \\ 
20126 Milano, ITALIA}
\email{{\tt mariagabriella.kuhn@unimib.it}}


\newtheorem{theorem}{Theorem}[section]

\newtheorem{proposition}[theorem]{Proposition}

\theoremstyle{definition}
\newtheorem{definition}[theorem]{Definition}
\newtheorem{remark}[theorem]{Remark}

\numberwithin{equation}{section}

\begin{document}


\begin{abstract}
Let $G=PGL(2,\Bbb Q_q)$. In  this paper
we shall investigate the group $\mcG$ of measurable
currents taking values in $G$.
The key observation is that $G$ is acting by automorphisms on a homogeneous
tree, which will play the role of the upper half plane in the case of
$PSL(2,\Bbb R)$.
Following the ideas of I.M. Gelfand, M.I. Graev and A.M. Vershik we shall 
construct an irreducible family of
representations of $\mcG$.
The existence of such representations 
depends deeply from the non-vanishing of the first cohomology group
$H^1(G,\pi)$ for a suitable infinite dimensional $\pi$.
\end{abstract}

\maketitle

\section{Introduction} 

Let $G$ be a locally compact group, $X$ any compact space and $\mcG^0$ the space
of all locally constant measurable functions $f:X\to G$.
The group structure of $G$
extends in a natural way to $\mcG^0$: for $f,g\in \mcG^0$ we let
$f\cdot g(x)=f(x)g(x)$.
Consider first $\mcG^0$ as a direct limit
of groups isomophic with $\underbrace{G\times \dots \times G}_{n\;\; times}$
and give it the natural topology coming from this structure. 
There are several approaches to the construction  
of continuous unitary irreducible
representations of $\mcG^0$ when $G$ is a Lie group. One method essentially 
 embeds $\mcG^0$ into the motion group
of a Hilbert space and uses the canonical projective representation
of the latter in the Fock space (see the papers of Araki \cite{2}, Guichardet
\cite{20}\cite{21}, Parthasarathy and Schmidt \cite{33}\cite{34} and Streater \cite{38}).

Our approach will follow the pioneering work of  I.M. Gelfand, M.I. Graev and A.M. Vershik \cite{18} and it is based on the existence of a semigroup
of positive definite functions called {\it canonical states}.

A semigroup of canonical states is a family of positive definite
functions $\Psi^\lambda(g)$ of the form
\begin{equation}
\Psi^\lambda(g)=\exp(\lambda(\psi(g)))\qquad\lambda>0
\end{equation}
where the infinitesimal generator $\psi(g)$ is a 
conditionally positive definite (briefly a c.p.d.)  function on $G$.
More references on c.p.d. can be found in the papers of P. Delorme
\cite{13} and \cite{14} and A. Guichardet \cite{21}. 

The existence of such a generator depends on the non-vanishing of the 
first cohomology group $H^1(G,\pi)$ where $\pi$ is an irreducible
representation of $G$ that cannot be separated
from the identity in the Fell topology (see \cite{24}). In this case one has 
$\psi(g)=-\frac12\norm{\beta(g)}^2$ for a suitable cocycle
$\beta\in H^1(G,\pi)$.

In \cite{18} and \cite{19}
 I.M. Gelfand, M.I. Graev and A.M. Vershik  
considered the case of  $PSL(2,\Bbb R)$, 
$SO(n,1)$ and $SU(n,1)$. 
In the same papers it was also showed that the representations constructed from
the semigroup $\Psi^\lambda$ are equivalent to those described
in the Fock model.

In this paper we shall consider the case of a  $p$-adic group of Lie type.

Let $\mcG^0$ be the group of locally constant functions  $f:X\to PGL(2,\Bbb Q_q)$
where $(X,\mu)$ is any compact measure space
(for technical reasons we shall assume
$\mu(X)<\frac12$). 

We shall consider two different ways to construct irreducible
representations of $\mcG^0$, both based on the existence of a non-trivial
cocycle $\beta\in H^1(G,\pi)$.

Denote by $(H_z,\pi_z)$ a complementary series representation of $\PG$
($z\in (0,\frac12)$). One way is based on the possibility of embedding
$(H_{z_1+z_2},\pi_{z_1+z_2})$ into $(H_{z_1}\otimes H_{z_2},\pi_{z_1}\otimes\pi_{z_2})$.
The other, apparently much more  closely related to $\beta$, is the construction of a representation $\Xi$ in the  Fock model. In our case $\beta$ will be an element
of $H^1(\PG,\St)$, where $\St$ is the Steinberg representation of $\PG$.

Theorem 6.3  investigates in detail the tensor product of two representations of the 
spherical series of $\PG$. Proposition 7.2 and 7.4 will allow us to construct a representation $\Pi$ of $\mcG^0$ as a direct limit of tensor products of representations of the complementary series of $\PG$.
Finally, in Theorem 7.5 we shall prove that
\begin{itemize}
\item $\Pi$ is irreducible
\item $\Pi$ is equivalent to $\Xi$.
\end{itemize}
In the last Section we shall show how to extend $\Pi$ to the group
$\mcG$ of {\it measurable bounded currents}.

The main difference between our case and those discussed in \cite{18}
concerns the semigroup of positive definite functions.
In the case of $PSL(2,\Bbb R)$ one has $\Psi^\lambda(g)=
\cosh^{-\lambda}(d_H(i,g\cdot i)/2)$ where $d_H$ is the hyperbolic distance
in the upper half plane between the two points $i$ and $g\cdot i$.

In the case of $PGL(2,\Bbb Q_q)$
one should expect that $\Psi^\lambda(g)=q^{-\lambda\dge}$ where
$\dge$ is tree distance between a choosen point $o$ and $g\cdot o$, however
the cocycle that corresponds to the c.p.d. function $-\dge$ is not {\it pure},
that is, it does
not take values in an irreducible
representation of $G$ (see \cite{28})  as in the case of $PSL(2,\Bbb R)$ 
and the semigroup $\Psi^\lambda$
is different from what expected.
The main analogy is about tensor products of
spherical representations (see Theorem 6.3 compared with \cite{36} and
Theorem II of \cite{35}).

\section{The tree associated to  $PGL(2,F)$}

Let $F$ be a commutative
non-archimedean local field. Let $\ord:F\to\Bbb Z\cup\{\infty\}$ be the
valuation on~$F$. Let $\mcO=\{x\in F:\ord(x)\ge0\}$ be the valuation
ring of~$F$, and let $\varpi\in\mcO$ be an element of valuation~1. Let
$\mcO^\times=\{x\in\mcO:\ord(x)=0\}$ denote the group of invertible
elements of the ring~$\mcO$. Let $q$ be the order of the residual
field~$\mcO/\varpi\mcO$, which equals~$p^r$ for some prime~$p$ and
some integer~$r\ge1$. Let $A\subset\mcO$ be a set of $q$ elements, one
of them~0, such that the canonical map $\mcO\to\mcO/\varpi\mcO$,
restricted to~$A$, is a bijection. Each element of~$\mcO$ is
expressible uniquely as the sum of a series
$a_0+a_1\varpi+a_2\varpi^2+\cdots$, where each $a_i$ is in~$A$.

Recall the construction of the Bruhat-Tits tree $\mcT$ associated
with~$H=GL(2,F)$ (see \cite{37}, p.~69; or \cite{16}, p.~127).
 Let $V=F^2$ denote the space
of all column vectors of length~2 with entries in~$F$. A lattice
in~$V$ is a subset of~$V$ of the form
$\{t_1v_1+t_2v_2:t_1,t_2\in\mcO\}$, where $\{v_1,v_2\}$ is a basis
of~$V$ over~$F$. If $\{v_1,v_2\}$ is the usual basis of~$V$, then the
corresponding lattice is~$\mcO^2$, and is denoted~$L_0$. If $L$ is a
lattice and if $g\in H$, then $g(L)$ is a lattice, and so
$H$ acts on the set of lattices. This action is clearly
transitive, and the stabilizer of~$L_0$ is the group $K=GL(2,\mcO)$ of
matrices with entries in~$\mcO$ and having determinant
in~$\mcO^\times$.  Two lattices $L,L'$ are called equivalent if
$L'=\lambda L$ for some $\lambda\in F^\times$. Let $[L]$ denote the
equivalence class of the lattice~$L$.  The Bruhat-Tits tree~$\mcT$ has
as vertex set the set of equivalence classes of lattices. Two distinct
lattice classes $[L]$ and~$[L']$ are adjacent if representative
lattices $L$ and~$L'$ can be found such that $\varpi L\subsetneqq
L'\subsetneqq L$. The tree~$\mcT$ is homogeneous of degree~$q+1$. 
Let $\Aut$ be the group of all automorphisms of $\mcT$.

The topology of pointwise convergence turns $\Aut$ into  a 
locally compact totally
disconnected topological group. Denote by
$\mcV$ the set of vertices of $\mcT$.
Let $\mcF$ be any finite subtree with vertex set $v_1\dots v_N$,
the sets $V_{\mcF}$  
\begin{equation}\label{intorni}
V_\mcF =\{g\in\Aut\;:g\cdot v_i =v_i \qquad\mbox{for all $v_i\in \mcF$}\}\period
\end{equation}
constitute, as $\mcF$ varies among all finite subtrees
of $\mcT$, a basis of neighbourhoods of the group identity $e$.

The above action of $H$ on~$\mcT$ gives a homomorphism $H\to\Aut$
with kernel $Z=\{\lambda I:\lambda\in F^\times\}$. Hence the quotient
group $PGL(2,F)=H/Z$ is isomorphic with a closed subgroup of $\Aut$
that we shall denote by 
$G$. We identify $G$ and~$PGL(2,F)$ throughout. The representations
of~$G$ correspond to, and are here frequently identified with,
representations of~$H$ which are trivial on~$Z$. 

From this point on we shall assume that {\it $G=PGL(2,\Bbb Q_q)$ where
  $q$ is a prime number different from $2$.}

It is natural 
to ask how the irreducible unitary representations~$\pi$ of~$\Aut$ behave when 
restricted to~$G$. When $\pi$ is spherical or special, the restriction
is known to remain irreducible (see \cite{15} and \cite{16}, p.~117), but the same is not true for other square integrable representations ( see \cite{8}).

\section{Representations of $G$ and Representations of $\Aut$}

Let $(H_\pi,\pi)$ be a conutinuous unitary representation of $G$.
A 
a $1$-cocycle $\beta:G\to H_\pi$
is a continuous function on $G$
satisfying the following
$$
\beta(g_1g_2)=\beta(g_1)+\pi(g_1)\beta(g_2)\;. 
$$
A  cocycle is trivial if it is cohomologous to zero, that is, 
$
\beta(g)=\pi(g)v-v\;
$
{for some $v\in H_\pi$}.

The cohomology group $H^1(G,\pi)$ is the quotient of the
vector space of all cocycles 
and the subspace of trivial cocycles.

The aim of this section is to identify  those irreducible representations
$(\pi, H_\pi)$ for which $H^1(G,\pi)$ is nontrivial.

Remember that $\pi$ is called {\it infinitesimally small} or a {\it core} 
if $\pi$ cannot be separated from the identity in the Fell topology.

In  1982 
Karpushev and Vershik \cite{24} proved the following
\begin{theorem}
Assume that $\beta:G\to H_\pi$
is a nontrivial cocycle and that
$(\pi,H_\pi)$ is irreducible. Then 
 $\pi$ cannot be
separated from the identity  in the Fell topology.
\end{theorem}

The  continuous cohomology
for p-adic Lie groups is known by results of Casselman \cite{11}.

Restricting our attention to the case of $G=PGL(2,\Bbb Q_q)$, we can  use
the description of
the unitary dual $\hat G$ given in \cite{17} (see also
 \cite{6} for more details) to search for those representations which are
infinitesimally small. It turns out that
the only possible  $\pi$ for which $H^1(G,\pi)$ is
nontrivial is $\St$,  the Steinberg representation, also called by many 
authors {\it special representation}.

We shall give here the description that is more convenient for our purposes.

\subsection{The boundary $\Omega$ and the special representation of $\Aut$.}

Let us denote by $d(v,v')$ the usual tree distance between vertices 
$v$ and $v'$.
A geodesic ray, or semi-infinite geodesic,
is a sequence of vertices $(v_0, v_1,\dots v_n,\dots)$
such that $d(v_i,v_{i+1})=1$ and $v_{i+2}\neq v_i$ for each $i$. 
The {\it boundary} $\Om$ is the set of equivalence classes of semi-infinite 
geodesics, two geodesics beeing equivalent if they coincide up to a shift.
For any given $v\in\mcT$ and $\om\in\Om$ there exists a unique geodesic
ray $[v,\om)$ representing $\om$ and starting at $v$.

  If $[v_0,\om)$ and $[w_0,\om)$ are two geodesic rays in the same class,
      the two semi-infinite geodesics coincide starting from a point called
{\it the confluent} and denoted by $v_0\wedge w_0$.

Fix once and for all the vertex $o$ corresponding to the lattice
$L_0=\mcO\oplus\mcO=\mcO^2$
and let $\mcK=\{g\in \Aut\st g\cdot o=o\}$.
Then $K=\mcK\cap G$ is a closed subgroup of $\mcK$ that  also acts
transitively on $\Om$.
For $g\in G$ and $\om\in\Om$ define the Busemann function
$$
h(o,g\cdot o; \om)= d(g\cdot o,x)-d(o,x) \quad 
$$
where $x=g\cdot o\wedge\om$ is the confluent of the geodesic rays $[o,\om)$
  and $[g\cdot o, \om)$.

Give $\Omega$ the natural topology as a subspace of the power space
$Map(\Bbb N, \mcV)$.
 This makes $\Omega$ compact
and totally disconnected, homeomorphic to the Cantor set.

Denote by $\nu$ the unique $\mcK$-invariant probability measure on $\Omega$ and 
by $\mcC^\infty(\Om)$ the space of complex locally constant functions.
The action of $G$ on $\mcT$ extends in an obvious way to an action on $\Om$.
The measure $\nu$ is quasi-invariant with respect to this action and
Radon-Nikodym derivative is given by
$$
\frac{d\nu(\inv g\om)}{d\nu\om}= q^{b(g,\om)}
$$
where
\begin{equation}\label{buseman}
b(g,\om)=-h(o,g\cdot o ;\om)
\end{equation}
and $b(g,\om)$ satisfies
\begin{equation}
b(g_1g_2,\om)=b(g_1,\om)+b(g_2,\inv{g_1}\om)\;.
\end{equation}

Fix $s\in\Bbb C$ and define a representation of $G$ on $\mcC^\infty(\Om)$ by
\begin{equation}
  \pi_s(g)F(\om)=F(\inv g\om)\left(\frac{s}{\sqrt q}\right)^{h(g\cdot o,o;\om)}
  =F(\inv g\om)\left(\frac{s}{\sqrt q}\right)^{-b(g,\om)}\;.
\end{equation}

The factor $\sqrt q$ gives the right normalization so that, when $|s|=1$,
$\pi_s$ is unitary with respect to the inner product
$\langle F_1,F_2\rangle=\int_\Om F_1(\om)\overline{F_2(\om)}d\nu$
and one obtains the {\it principal series} acting on
$L^2(\Om,d\nu)$.
 When
$s$ is real with $1/\sqrt{q}<|s|<\sqrt{q}$ the representations $\pi_s$ are
unitarizable 
 and the
{\it complementary} spherical series act on $H_s$, the completion of $\Kk$
with respect to another  suitable inner product.
Principal and complementary series make up the {\it spherical series}: they
are the only irreducible representations admiting a nonzero $\mcK$-invariant 
vector. Since $(\Aut, \mcK)$ is a Gelfand pair, the subspace of $\mcK$-invariant
vectors in $\pi_s$ is one dimentional.

The following Theorem is known from the work of P. Cartier \cite{7}, where
spherical functions of $G$ are computed in terms of the action on $\mcT$.
We shall give here a proof that will be needed later to deal with the Steinberg
representation.
\begin{theorem} The principal or complementary series
of $\Aut$ restrict to the  {\it spherical} principal or complementary series
of $G$.
\end{theorem}
\begin{proof}
Let $\chi_s: F^\times\to {\Bbb C}^\times$ be the quasi-character
$a\mapsto s^{\ord(a)}$ of $F^\times$. 

Then it is routine to see that 
$\pi_s$ 
is the principal series representation
$\rho_s=B(\chi_s,\chi_{s^{-1}})$ as defined in Bump's book~\cite{6}( p.471).
 Indeed,
let $\omega_0$ be the class of the geodesic $(g_0o,g_1o,\ldots)$,
where $g_n=\begin{pmatrix}
1&0\\
0&\varpi^n\end{pmatrix}$
 for $n\in\Bbb N$. The set of 
$g\in G$ such that
$g\omega_0=\omega_0$ is the subgroup~$P=\left(
\begin{matrix}
a & b\\
0 & c
\end{matrix}\right)
$ of upper-triangular matrices in~$G$.
We define $T:\mcC^\infty(\Omega)\to V_s$, the representation space 
of~$\rho_s$ by
$$
(TF)(g)=F(g^{-1}\omega_0)\Bigl({s\over\sqrt{q}}\Bigr)^{h(o,g\cdot o;\omega_0)}.
$$
With the notation of I.M.Gel'fand, M.I.Graev, I.I. Pyatetskii-Shapiro
 \cite{17} $\pi_s$ is  the representation  induced from
 the character $\chi^s$ of $P$, obtained from $\chi_s$ by letting
$\chi^s
\left(\begin{matrix}
a & b\\
0 & c
\end{matrix}\right)
=\chi_s(a)/\chi_s(b)$.
\end{proof}
\begin{remark}
It is important to observe that our character $\chi_s$ of $F^\times$
  {\it is trivial}
on $\mcO$: when $\chi_s$ is not trivial on $\mcO$ one can still obtain
irreducible representations of $G$ which {\it are not}  restrictions
to $G$ of irreducible representations of $\Aut$.
It can be proved (see 
\cite{9}and \cite{16}) that some {\it discrete series} 
representations of $\Aut$ decompose, when restricted to $G$, as a direct 
integral of representations
 coming from unitary induction
of characters that are not trivial on $\mcO$.
In particular the {\it discrete series} of $\Aut$ {\it do not restrict}
in general
to the discrete series of $G$. We shall see that a remarkable exception
is one of the two {\it special representations}.
\end{remark}

\subsection{The special representation $\sp$ of $\Aut$}

This representation of $G$ appears in the limit, as $s\to\sqrt{q}$ or 
$s\to1/\sqrt{q}$ of the complementary spherical series.
More precisely when $s$ approaches $s_1=1/\sqrt{q}$ the representation $\pi_s$
approaches a representation $\pi_{s_1}$ which  has the irreducible representation
$\sp$ as a subrepresentation and the trivial representation  as a quotient. When
$s$ approaches $s_0=\sqrt{q}$ the limit representation $\pi_{s_0}$ contains
the trivial representation and the special representation $\sp$ is realized
on the quotient space. 
\begin{theorem}
The Steinberg representation $\St$ of $G$ is the restriction to $G$
of the special representation $\sp$ of $\Aut$.
\end{theorem}
\begin{proof}
Since the spherical series of $G$ can be described only by means of the action
of $G$ on the boundary $\Om$ the same is true for their limits. 
Since  the Steinberg representation of $G$ can be 
obtained taking the limits of the spherical complementary series, it can also
be  regarded as the limit of the restriction
to $G$ of the complementary series of  $\Aut$. 
In this last case the representations approach $\sp$,  the
special representation of $\Aut$  (see \cite{32}). 

The representation $\pi_{s_1}$ of $\Aut$ is the  quasi-regular representation
 ( it is induced from the trivial representaion
of the stabilizer of $\omega_0$).
Hence the restriction of $\pi_{s_1}$ to $G$ is the quasi-regular representation
of $G$ induced from the trivial representation of $P$.
 They both act on $\CC$ according to the rule
\begin{equation}\label{quasireg}
\pi_{s_1}(g)F(\om)=q^{b(g,\om)}F(\inv g\om)\;.
\end{equation}
Since $q^{b(g,\om)}=\frac{d\nu(\inv g\om)}{d\nu(\om)}$ 
is the Radon-Nikodym derivative of the action,  $\pi_{s_1}$
preserves $\int_\Om F(\om)d\nu(\om)$.
Hence the subspace $\mcC_0$
 consisting of functions having zero integral is invariant
by $\pi_{s_1}$. The scalar product $\langle F, F\rangle_{s_1}=\int_{\Om\times\Om}
F(\om_1)F(\om_2)d\nu(\om_1)d\nu(\om_2)$ is also 
preserved by $\pi_{s_1}$ and the subspace
consisting of all functions $F$ such that $\langle F,F\rangle_{s_1}=0$ is $\mcC_0$.
In this context the quotient representation $\CC/{\mcC_0}$ is the trivial one 
and the Steinberg representation of $G$ is the restriction to $G$ of
the special representation of $\Aut$, acting on $\mcC_0$.

The other limit, when $s$ approaches $s_0=\sqrt{q}$, gives the same 
representation of $\Aut$ and hence the same restriction to $G$. We shall
fell free to use any of these two descriptions. In particular the one coming
from $\pi_{s_0}$ is more convenient for our purposes.

The representation $\pi_{s_0}$  of $\Aut$ 
is acting on $\CC$ according to the formula
\begin{equation}
\pi_{s_0}(g)F(\om)=F(\inv g\om)\;.
\end{equation}

Consider two distinct elements $\om_1=(o,v_1,v_2,\dots)$
and $\om_2=(o,w_1,w_2,\dots)$ of $\Om$ both represented by
two geodesic rays starting at $o$. If $v_1\neq w_1$ set
$|\om_1\wedge\om_2|=0$. If not, denote by $k+1$ the smallest integer such that
$v_{k+1}\neq w_{k+1}$ (so that $v_j=w_j$ for all $j$ with $0\leq j\leq k$)
and let $|\om_1\wedge\om_2|=k$.

Define the inner product
\begin{equation}\label{innerprod}
\langle F,F\rangle_{s_0}=\frac{q^2-1}{q}
\int_\Omega\int_\Omega|F(\omega)-F(\omega')|^2
q^{2|\omega\wedge\omega'|}
d\nu(\omega)d\nu(\omega')
\end{equation}
and consider the quotient space $\mcC_0=\CC/\{\langle F,F\rangle_{s_0}=0\}$.
It is known  that the measure 
$q^{2|\omega\wedge\omega'|}d\nu(\omega)d\nu(\omega')$ is invariant under the
diagonal action $g\mapsto(g\om_1, g\om_2)$ of $G$ on $\Om\times\Om$
(see \cite{23}). Hence the above inner product \eqref{innerprod}
is preserved by $\pi_{s_0}$. The special representation can be realized
as $\pi_{s_0}$ acting on $H_{s_0}$ the completion of $\mcC_0$ with respect to
$\langle F,F\rangle_{s_0}$
(see for example \cite{26} section 3).
\end{proof}

Let $v$ be any element of $H_{s_0}$.
From now on we shall denote, for simplicity,

\begin{align*}
&  \norm{v}^2=\langle v,v\rangle_{\sp}=\langle v,v\rangle_{s_0}=\\
&\frac{q^2-1}{q} \int_\Omega\int_\Omega|v(\omega)-v(\omega')|^2
q^{2|\omega\wedge\omega'|}d\nu(\omega)d\nu(\omega')
\end{align*}

and we shall identify
the Steinberg representation $\St$ with  
the restriction to $G$ of the special
representation $\sp$ of $\Aut$ realized on $H_{s_0}$.
\section{The cocycle}

The continuous 1-cohomolgy of $\Aut$ had been computed by C. Nebbia
\cite{30}. The same result was also obtained (but not published) by O. Amann
in his Master Thesis \cite{1}.

We shall now give the description that is more convenient for our purposes,
starting from the description of a nontrivial element of $H^1(G,\St)$.

For every $g\in G$, let $b(g)$ the function on $\Om$ defined by the rule
\begin{equation}
  b(g)(\om)=b(g,\om)\;.
  \end{equation}
Since $q^{b(g,\om)}$ is a Radon-Nikodym derivative (see \eqref{quasireg})
$b(g,\om)$ satisfies
\begin{equation}\label{cocycle}
b(g_1g_2,\om)=b(g_1,\om)+b(g_2,\inv{g_1}\om)\;.
\end{equation}

\begin{theorem}
  The assignment $g\mapsto b(g)$ defines an unbounded
  cocycle on $G$ with values in $\St$.
\end{theorem}
\begin{proof}
  Identify $\St$ with $\pi_{s_o}$ acting on $H_{s_0}$.
  Since $\pi_{s_0}(g)F(\om)=F(\inv g\om)$, $b(g)$ is a cocycle
  by \eqref{cocycle}.
By Theorem 5 of \cite{26} we have:
\begin{equation}
\frac12\norm{b(g)}^2=\dge +\frac{2q}{q^2-1}\left(q^{-\dge}-1\right)
\end{equation}
showing that $b(g)$ is clearly unbounded.
\end{proof}

It should be noticed that another quite useful
realization of the special representation, as the first $L^2$-cohomology of the tree, is given by A.Borel in \cite{5}.

{\bf The canonical semigroup}

At this point, in analogy with the case of $PSL(2,\RR)$,  we shall put
\begin{equation}
\Psi^\lambda(g)=\exp(-\frac12\lambda\|b(g)\|^2)\qquad\lambda>0
\end{equation}

Let, as usual, $\mcG^0$ denote the group of step functions $f:X\mapsto G$
endowed with the direct limit topology. Let $\mu$ be a positive Borel measure
on $X$ such that $\mu(X)=M<\frac12$
We modify a little bit
the norm on $H_0$ in order to avoid tedious changes of parameters for
the complementary series of $G$. For $v\in H_0$ set
\begin{equation}\label{hm}
\norm{v}^2_q=\log(q)\langle v,v\rangle\;.
\end{equation}
Denote by $H_q$ the Hilbert space obtained from $H_0$ modifying the norm
according to the formula above \eqref{hm}.
Let $H^X=\int^\oplus H_x d\mu\simeq L^2(X,d\mu)\otimes H_q$
 denote the direct integral of spaces
$H_x=H_{q}$ a.e. $[\mu]$.
 So that $H^X$ 
is the completion of locally constant
measurable mappings $\bold v:X\to H_{q}$ with respect to the norm
\begin{equation}\label{norma}
\norm{\bold v}^2=\log(q)\int_X\langle{v(x)},v(x)\rangle\; d\mu(x)\period
\end{equation}
Define a cocycle  $b^X :\mcG^0\to H^X$
and a representation $\St^X:\mcG^0\to H^X$ by the rule
\begin{align}
b^X(\xi)(x)&=b(\xi(x)) \\
\St^X(\xi)(\bold v)(x) &=\St(\xi(x))v(x) 
\end{align}

Let $\mcH=\Exp(H^X)$ be  the Fock space constructed from $H^X$ and
let $\Xi:\mcG^0\to\mcH$ defined as follows: 
\begin{equation}\label{ourexp}
\begin{aligned}
&\Xi(\xi)(\Exp(\bold v))=\\
&e^{\displaystyle{(-\frac12\norm{b^X(\xi)}^2-
\langle\St^X(\xi)v,b^X(\xi)\rangle)}} 
\cdot\Exp(\St^X(\xi)\bold v+b^X(\xi))\period
\end{aligned}
\end{equation}

Usually $\Xi$ is a unitary {\it projective} representation of $\mcG^0$, but
in our case, since $b(g)$ is real and $\St(g)$ preserves the subspace of real functions, $\Xi$ is a true (unitary) representation of $\mcG^0$.

\begin{theorem}
The representation $\Xi$ in the Fock space $\mcH$
is irreducible.
\end{theorem}
\begin{proof}
  Argue as in the proof of Theorem 4.2 of \cite{25} and apply
 an irreducibility criterium taken from Ismagilov~\cite{22} which works also in this case.
\end{proof}
\section{Another Realization of $\Xi$.}

Another possible realization of $\Xi$ is based
on the tensor product Theorem for spherical representations.
There are strong analogies between the case of $PSL(2,\Bbb R)$, 
$PGL(2,\Bbb Q_q)$ and $\Aut$.
The aim of this section is to exhibit  these analogies.

\subsection{The case of $SL(2,\Bbb Q_q)$}

Tensor products of irreducible representations of $SL(2,\Bbb Q_q)$ have been 
studied  by many authors. The case of two principal series representations
had been considered by R.P. Martin \cite{27} while many other cases
are considered
 in the papers of C. Asmuth and J. Repka (see \cite{3}, \cite{4}).

In particular we are interested in the case of tensor product of two
{\it complementary} series representations which is treated in 
Theorem 3.6 of \cite{3}. In terms of our 
parameters the result is the follwing: 
\begin{theorem}[C. Asmut, J. Repka]
Take two
complementary series of $SL(2,\Bbb Q_q)$, say $\pi_s$ and $\pi_{s'}$.
Assume for semplicity that $1<s,s'<\sqrt{q}$ and
consider the tensor product $\pi_{s,s'}=\pi_s\otimes\pi_{s'}$ of the two spherical
representations.
\begin{itemize}
\item
If $|ss'|\leq\sqrt{q}$ then
\begin{equation*}
\pi_{s,s'}\simeq \pi_1\otimes\pi_1
\end{equation*}
where $\pi_1$ is a principal series representation.
\item
If $|ss'|>\sqrt{q}$ then
\begin{equation*}
\pi_{s,s'}\simeq (\pi_1\otimes\pi_1)\oplus \pi_c
\end{equation*}
where $\pi_1$ is as before while $\pi_c$ is the complementary series
representation corresponding to
the parameter $c=\frac{ss'}{\sqrt{q}}$.
\end{itemize}
\end{theorem}

\section{The case of $PGL(2,\Bbb Q_q)$}
Let us turn now to our case. Besides some technical devices in passing from
$SL(2,\Bbb Q_q)$ to  $PGL(2,\Bbb Q_q)$, the crucial point is that, in order
to show that $\Xi$ is equivalent to the representation constructed by direct
 limits, we need to know exactly
how $\pi_c$ embeds into $\pi_{s,s'}$.
Our strategy is the following: we know exactly how this embedding goes for
$\Aut$, the group of automorphisms of $\mcT$. 
Hence we shall first decompose $\pi_s\otimes\pi_{s'}$ as 
{\it representation of $\Aut$}. Finally, using the results about 
the restrictions of the spherical series and the special representation, we 
shall derive the formulas needed for $G$.
 
\subsection{The case of $\Aut$}

We are interested here in the representations of $\Aut$ that appear in the
Plancherel formula but that are not contained in the regular representation.
We know that only the representations of the principal spherical series
have this property (see \cite{31}, \cite{32}, \cite{16}).
Denote by $\one$ the function identically $1$ on $\Om$.
Let $s$ be a nonzero complex number.
The {\it principal series representations}
correspond to the case $|s|=1$ and the
'endpoints' are obtained when  $s=\pm1$. We distinguish between 
$|s|=1$ and  $s\neq\pm1$ because in this last case the spherical functions
have a slightly different behavior. Namely,
when $s\neq\pm1$ one has

%
\begin{equation}\label{sferiche}
\varphi_n(s)=\langle \pi_s(g){\one},{\one}\rangle
=\q^{-n/2}\bigl(c(s)s^n+c(s^{-1})s^{-n}\bigr)
\;\text{\rm if\ }\dge=n 
\end{equation}
where
\begin{equation}\label{HCS}
c(s)=\frac{qs-\inv s}{(q+1)(s-\inv s)}
\end{equation}
is the Harish-Chandra $c$- function.
When $s=\pm1$ one has
\begin{equation}
\varphi_n(s)={s^n\over(\q+1)\q^{n/2}}\bigl((\q+1)+(\q-1)n\bigr)\;.
\end{equation}
We recall that the expression \eqref{sferiche} is also valid for the representations
of the {\it complementary series} whose 'endpoints' are obtained
for $|s|=\sqrt{q}$ and $|s|=1/\sqrt q$.

Denote by $ds$ the normalized Haar measure on $\Bbb T$. The continuous
part of the Plancherel measure for $\Aut$ is given by
\begin{equation}
d\mu_p(s)=\frac{q}{2(q+1)}\frac1{|c(s)|^2}ds\;.
\end{equation}

The following proposition is standard (see for example \cite{35}):
\begin{proposition}
Let $\pi_{s,s'}=\pi_s\otimes\pi_{s'}$ be the tensor product
of two representations of the spherical
series acting on $H_s\otimes H_{s'}$. 
Then
\begin{equation}
H_s\otimes H_{s'}=H_K\oplus H_N
\end{equation}
and $\pi_{s,s'}$ is the direct sum of two representations, 
$\pi_K$ and $\pi_N$, acting respectively on $H_K$ and $H_N$
 with the following properties:
 $H_N$ does not contain any nontrivial $K$-invariant
vector while $H_K$ decomposes by means of the spherical series only.
Furthermore $H_K$ equals the closure
in~$H_{s}\otimes H_{s'}$
 of the linear span of~$\{\pi_s\otimes\pi_{s'}(g)(\one\otimes\one):g\in G\}$.
\end{proposition}

The decomposition of $H_N$ is described in \cite{8}. The only thing
that we need to know here is that $H_N$ is the direct sum of
the discrete series of $\Aut$.

Decomposing $H_K$ is the same as decomposing the positive definite
function $\varphi_s(g)\varphi_{s'}(g)$ into the sum, or direct integral, of
spherical functions. The following Theorem is taken from \cite{10}

\begin{theorem}Assume that 
 $s$ and~$s'$ are
parameters corresponding to the principal or complementary series, with
$s,s'\not\in\{\pm\sqrt{q},\pm1/\sqrt{q}\}$. Write, as usual $n=\dge$. Then 
\begin{itemize}
\item If $|ss'|\leq\sqrt{q}$ 
\begin{equation}
\varphi_{s}(g)\varphi_{s'}(g)=\int_{\Bbb T}K(s,s',t)\varphi_t(g)d\mu_p(t)\;.
\end{equation}
\item If $|ss'|>\sqrt{q}$
\begin{equation}\label{prspfct}
\varphi_{s}(g)\varphi_{s'}(g)=A\varphi_c(g)+\int_{\Bbb T}K(s,s',t)
\varphi_t(g)d\mu_p(t)\;.
\end{equation}
where $c= ss'/\sqrt{q}$, $A$ is a positive number less than $1$
given by the Harish-Chandra $c$-function:
\begin{equation}
A=c(s)c(s')/c(ss'/\sqrt{q})
\end{equation}
and $K$ is the kernel
\begin{equation}
K(s_1,s_2,s_3)=
{\q(\q-1)\over(\q+1)^2}
{\prod_{j=1}^3\bigl\{\bigl(1-{s_j^2\over \q}\bigr)
\bigl(1-{s_j^{-2}\over \q}\bigr)\big\}
\over
\prod_\epsilon(1-{1\over\sqrt{\q}}
s_1^{\epsilon_1}s_2^{\epsilon_2}s_3^{\epsilon_3})}\;.
\end{equation}
\end{itemize}
The above product is over all eight
$\epsilon=(\epsilon_1,\epsilon_2,\epsilon_3)\in\{+1,-1\}^3$ and,
 in both cases, $K(s,s',t)>0$
for all $t\in\Bbb T$.\par
\end{theorem}

Let us turn now to the case of our $G=PGL(2,\Bbb Q_q)$:
\begin{theorem}
Assume that $\pi_s$ and $\pi_{s'}$ are representations from the spherical
series of $G$ and form the tensor product $\pi_{s,s'}=\pi_s\otimes\pi_{s'}$.
\begin{itemize}
\item
  If $|ss'|\leq\sqrt{q}$ the representation $\pi_{s,s'}$ is
  contained  in the direct sum of copies of the
  the regular representation of $G$
\item
If $|ss'|>\sqrt{q}$ 
the representation $\pi_{s,s'}$ splits into the direct sum of 
two pieces: one 
is the complementary series corresponding to the parameter $ss'/\sqrt{q}$, the
other is contained in the direct sum of copies of the
regular representation.
\end{itemize}
\end{theorem}
\begin{proof}
  Since taking the restrictions of the spherical
  principal/complementary series of $\Aut$
  we get  the spherical
  principal/complementary series of $G$, the only thing that we 
need to be prove is that the restriction to $G$ of a square integrable
representation of $\Aut$ is contained in the regular representation of $G$.
This is certainly true if the representation of $\Aut$ is a direct integral
of the spherical principal series.
When the representation belongs to the discrete series of $\Aut$ it may happen
that some { discrete series} of $\Aut$ decompose, when restricted to $G$, as a direct 
integral of representations
 coming from unitary induction
of characters that are not trivial on $\mcO$ (see 
\cite{9}and \cite{16}). 
Nonetheless the restriction to $G$ of a discrete series representation of $\Aut$
is still contained in the direct sum of a finite number of copies of the regular representation of $G$: see Proposition 2.2
of \cite{9}.
\end{proof}

From this point   {\it principal series} will stay for both
spherical principal series (representations unitarly
induced from characters trivial
on $\mcO$) and representations unitarly
induced from characters non trivial on $\mcO$.

We summarize with the following:
\begin{theorem}\label{embedding}
Let $\pi_{s,s'}=\pi_s\otimes\pi_{s'}$ be the tensor product of two complementary
series representations of $G$ acting on $H_{s,s'}$. Assume that
$|ss'|>\sqrt q$ and let $c=\frac{ss'}{\sqrt q}$. Then the representation
$\pi_{s,s'}$ contains the complementary series representation $\pi_c$ with 
multiplicity one.
Moreover the vector $\one\otimes\one$ can be written as the orthogonal sum
of two vectors
$$
\one\otimes\one=v_c\oplus w
$$
in such a way that
$$
\langle\pi_s\otimes\pi_{s'}(g)v_c,v_c\rangle=
          {\frac{c(s)c(s')}{c(ss'/\sqrt q)}}\varphi_{ss'/\sqrt q}(g)
$$          

and the closed linear span of $\pi_s\otimes\pi_{s'}(g)w$ is contained in the regular representation of $G$.
Observe that the map $\one\to\sqrt{\frac{c(ss'/\sqrt q)}{c(s)c(s')}}
v_c$ extends to an isometric embedding of $H_c$ into $H_{s}\otimes H_{s'}$     
\end{theorem}

\begin{remark} Since $(G,K)$ is a Gelfand pair we know that the vector $v_c$
  is uniquely determined up to scalars of absolute value one.
\end{remark}

\section{The representation of $\mcG^0$}

\begin{definition}
  Let $A$ be any Borel subset of $X$ and let $G_A$ be the group of functions
  $f:X\to G$ that are constant on $A$ and identically one on the complement of $A$. Obviously one has $G_A \simeq G$.
  A  partition $\rho$ of $X=\cup_{i=1}^n A_i$
  into finitely many disjoint Borel subsets is called admissible.
  Given two admissible partitions $\rho_1$ and $\rho_2$ we say that
  $\rho_1<\rho_2$ if $\rho_2$ is a refinement of  $\rho_1$.
  For any admissible partition $\rho$  of $X=\cup_{i=1}^n A_i$,
  $G_\rho$ will denote the group of functions $g:X\to G$ that are constant on each of the subsets
  $A_i$. Hence $G_\rho\simeq G_{A_1}\times\dots\times G_{A_n}$.
  For $\rho_1<\rho_2$ there is a natural embedding of $G_{\rho_1}$
  into $G_{\rho_2}$, so that we may define the inductive limit
  $$
  \mcG^0=\lim_{\rightarrow}G_\rho
  $$
\end{definition}  

\subsection{Representations of $\mcG^0$}
It is convenient at this point to change the parametrization of the spherical series, namely we set
$$s_i=q^{\frac12-z_i}$$
so that $\pi_{s_i}$ will be denoted by $\pi_{z_i}$. We are interested in the case
$0<z_i<\frac12$. With this new parametrization
the multiplication formula \eqref{prspfct}
reads as
\begin{equation}\label{newmult}
  \varphi_{z}(g)\varphi_{z'}(g)=\frac{c(z)c(z')}{c(z+z')}\varphi_{z+z'}(g)+
  \int_{\Bbb T}K(z,z',t)
\varphi_t(g)d\mu_p(t)
\end{equation}
where the Harish-Chandra $c$ function
becomes
\begin{equation}\label{HaCh}
c(z)=\frac1{(q+1)}\frac{q^{1-z}-q^{z-1}}{q^{-z}-q^{z-1}}\;.
\end{equation}

With respect to this parameters {\it the product $ss'/\sqrt q$ is turned into $z+z'$,}
which is more convenient for our purposes.

Let $\mu$ be a positive Borel measure on $X$ such that $\mu(X)=M<\frac12$.
We shall now associate to any admissible partition of $X=\cup_{i=1}^n A_i$ 
an irreducible representation
of $G_\rho$. 

Denote by $H_0$ the one-dimentional space on which the identity representation
of the trivial group acts.
Let 
$\rho\;:X=\cup_{i=1}^n A_i$ be any admissible partition. Set $z_i=\mu(A_i)$.
 Let $\pi_{z_i}$ be the complementary series
 representation corresponding to this parameter
 acting on $H_{z_i}$.
 
 Define $H_\rho=H_{z_1}\otimes\dots H_{z_n}$ and
 an irreducible
representation $\pi_\rho$ of $G_\rho$ by the rule 
$\pi_\rho(g_1,\dots g_n)=\pi_{z_1}(g_1)\otimes\dots\otimes\pi_{z_n}(g_n)$ .

\begin{proposition}\label{morfismi}
  Let $A$ be the disjoint union of measurable sets
 $A_1\dots A_n$ with
$\mu(A)=z$ and $\mu(A_i)=z_i$. Let  
$\pi_z$, respectively $\pi_{z_i}$, denote the complementary
series representations of $G$ acting on $H_z$, respectively on $H_{z_i}$.
The representation $\pi_{z_1}\otimes\dots\otimes\pi_{z_n}$
of $G$ splits into of  
the orthogonal sum of two pieces:
\begin{equation}
\pi_{z_1}\otimes\dots\otimes \pi_{z_n}=\pi_z \oplus \pi'
\end{equation}
where $\pi_z$ is the complementary series representation corresponding to the 
parameter $z$ and 
$\pi'$ decomposes by means of the principal and discrete series only.
In particular there exists an isometric embedding of
$H_z$ into $H_{z_1}\otimes \dots \otimes H_{z_n}$
that commutes with the action of $G$. 
\end{proposition}
\begin{proof}
  Theorem \eqref{embedding}
says that the statement is true for $n=2$ 
Assume that $z=z_1+z_2+z_3$, multiply both sides of \eqref{newmult}
by $\varphi_{z_3}(g)$ and apply again Theorem \eqref{embedding}:
\begin{equation}\label{pprspfct}
\varphi_{z_1}(g)\varphi_{z_2}(g)\varphi_{z_3}(g)=\frac{c(z_1)c(z_2)}{c(z_1+z_2)}
\frac{c(z_1+z_2)c(z_3)}{c(z_1+z_2+z_3)}\varphi_{z_1+z_2+z_3}(g)+
\lambda(g)
\end{equation}
Where 
\begin{multline}
\lambda(g)=\varphi_{z_3}(g) \int_{J}K(z_1,z_2,t)\varphi_{\frac12+it}(g)d\mu(t)\;+\\
\frac{c(z_1)c(z_2)}{c(z_1+z_2)}
\int_{J}K(z_1+z_2,z_3,t)\varphi_{\frac12+it}(g) d\mu(t)
\end{multline}
Since the representations
that are weakly contained in the regular representation 
are characterized by the decay of their matrix coefficents 
(see \cite{12}),
 the tensor product of a uniformly bounded representation and a 
representation weakly contained in the regular is still weakly
contained in the regular and we may conclude that no complementary
series appears in the decomposition of $\lambda$.
In particular there exists a vector  $v_{z_1,z_2,z_3}$
on $\Om\times\Om\times\Om$ 
such that the map
\begin{equation}
\one_{z_1+z_2+z_3}\to\sqrt{\frac{c(z_1+z_2+z_3)}{c(z_1)c(z_2)c(z_3)}}v_{z_1,z_2,z_3}
\end{equation}
extends to an isometric embedding $j_{z_1,z_2,z_3,\, z}$
of $H_z$ into
 $H_{z_1}\otimes
H_{z_2}\otimes H_{z_3}$.
Repeated applications of the above arguments conclude the proof.
\end{proof}
\begin{remark}
  For any $z$ with $0<z<\frac12$ denote by $\one_z$ the unique
  positive $K$-invariant function of norm one in $H_z$.
  The above Theorem shows that, for any partition
  $\rho$ such that $z=z_1+z_2\dots+z_n$
  it is possible to choose vectors $v_\rho=v_{z_1,\dots,z_n}\in H_{z_1}\otimes \dots \otimes H_{z_n}$ such that
  \begin{equation}\label{amen}
    \one_{z_1}\otimes \dots \otimes \one_{z_n}=v_{\rho}\oplus w
  \end{equation}
  where the above sum is an orthogonal sum and the closed linear span
  of $\pi_{z_1}\otimes \dots \otimes \pi_{z_n}(g)w$ does not contain any
  representation of the
  complementary series. Moreover one has
  $$
  \langle \one_{z_1}\otimes \dots \otimes \one_{z_n}(g)v_{\rho},v_\rho\rangle=
          \frac{c(z_1)\dots c(z_n)}{c(z_1+\dots+z_n)}
          \varphi_z(g)
  $$
  The vectors $v_\rho$ are uniquely determined up to scalars of absolute value one and  the assignment $\one_z\to \sqrt{\frac{c(z_1+\dots+z_n)}{c(z_1)\dots c(z_n)}} v_\rho$ defines an isometric embedding of $H_z$ into $\mcH_\rho$.
\end{remark}
\begin{proposition}
  Given three admissible partitions
  $\rho_0<\rho_1<\rho_2$ there exist isometric $G$ embeddings $j_{\rho_1,\rho_0}$,
  $j_{\rho_2,\rho_1}$ such that the following diagram commutes:
\begin{equation*}
\begin{tikzcd}[row sep=tiny]
  & H_{\rho_1}\arrow[dd,"j_{\rho_2,\rho_1}"] \\
    H_{\rho_0}\arrow[ru,"j_{\rho_1,\rho_0}" pos=0.7] \arrow[rd,"j_{\rho_2,\rho_0}"' pos=0.7] \\
  & H_{\rho_2}
\end{tikzcd}
\end{equation*}


  \end{proposition}
\begin{proof}


  Assume first that $\rho_0=\{A\}$. Split $A$ into two pieces, say $A_1$ and $A_2$ and finally split $A_1$ into two more pieces, say $A_{1,1}$ and $A_{1,2}$.
  Let $\rho_1$ and $\rho_2$ be the corresponding partitions and let
  $v_{\rho_i}$ as in \eqref{amen}.
Write $\mu(A)=\mu(A_1)+\mu(A_2)=z=z_1+z_2$.
Say that $z_1=\mu(A_{1,1})+\mu(A_{1,2})=w_1+w_2$.
  Write, as
  in \eqref{amen}
  \begin{align}
    &\one_{w_1}\otimes\one_{w_2}\otimes\one_{z_3}=v_{\rho_2}\oplus w\\
    &\one_{z_1}\otimes\one_{z_2}=v_{\rho_1}\oplus \tilde w\\
    &\one_{w_1}\otimes\one_{w_2}= v_{1,2} \oplus \bar w
  \end{align}
  where, in the last equation, $v_{1,2}$ is chosen in such a way that
  the map $\one_{z_1}\to \sqrt{\frac{c(z_1)}{c(w_1)c(w_2)}} v_{1,2}$ defines an isometric embedding
  of $H_{z_1}$ into $H_{w_1}\otimes H_{w_2}$. Hence
  
  $$
  \one_{w_1}\otimes\one_{w_2}\otimes\one_{z_2}=v_{1,2}\otimes\one_{z_2}\oplus W
  $$
  where the closed linear span of
  $\pi_{w_1}\otimes\pi_{w_2}\otimes \pi_{z_2}(g)W$ does not contain any
  representation of the
  complementary series.
  Compute
  $$
  \langle  \pi_{w_1}\otimes\pi_{w_2}\otimes\pi_{z_2}(g)v_{\rho_2},v_{\rho_2}\rangle
  =\frac{c(w_1)c(w_2)c(z_2)}{c(w_1+w_2+z_2)}\varphi_z(g)
  $$
  while
  \begin{align}
  &\langle  \pi_{w_1}\otimes\pi_{w_2}\otimes\pi_{z_2}(g)v_{1,2}\otimes\one_{z_3},
    v_{1,2}\otimes\one_{z_2}\rangle=\\
    &\langle  \pi_{w_1}\otimes\pi_{w_2}(g)v_{1,2},v_{1,2}\rangle\cdot
    \langle\pi_{z_2}(g)\one_{z_2},\one_{z_2}\rangle=\\
    & \frac{c(w_1)c(w_2)}{c(z_1)}
    \varphi_{z_1}(g)\varphi_{z_2}(g)=\\
    & \frac{c(w_1)c(w_2)}{c(z_1)}\frac{c(z_1)c(z_2)}
        {c(z_1+z_2)}\varphi_z(g)=\\
   &=\frac{c(w_1)c(w_2)c(z_2)}{c(w_1+w_2+z_2)}\varphi_z(g)     
  \end{align}
  
  Repeated applications of the above argument conclude the proof.
    \end{proof}

\subsection{The representation space for $\mcG^0$}

For any admissible partition $\rho$ construct the Hilbert space $\mcH_\rho$.
Theorem \eqref{morfismi} above ensures that for any pair 
 $\rho_1<\rho_2$ there exist morphisms
 $j_{\rho_2,\rho_1}:\mcH_{\rho_1}\to\mcH_{\rho_2}$ which commute with the
$G_{\rho_i}$ action and also satisfy the compatibility condition
 $j_{\rho_3,\rho_2}\cdot j_{\rho_2,\rho_1}=j_{\rho_3,\rho_1}$
 for any $\rho_1<\rho_2< \rho_3$.
We can now define $\mcH^0$ to be the inductive limit of Hilbert spaces:
\begin{equation}
\mcH^0=\lim_{\rightarrow}\mcH_\rho\;.
\end{equation}
We recall that an element $v$ of $\mcH^0$ is an equivalence 
class of vectors $[v_\rho]$ with  $v_\rho\in\mcH_\rho$ and with the
 following property:
there exists $\rho_0$ (depending on $v$) such that for every 
admissible partition $\rho>\rho_0$ one has $v_\rho=j_{\rho,\rho_0}v_{\rho_0}$.
One can take 
\begin{equation*}
\Norm{v}_{\mcH^0}=\Norm{v}_{\mcH_{\rho_0}}\;.
\end{equation*}
Let $\mcH$ be the completion of $\mcH^0$ with respect to this norm.
Let now
$\xi\in\mcG^0$ and $v\in\mcH^0$. In order to define $\Pi(\xi)v$ observe 
that
 there exist an  admissible partition $\rho$ 
such that $\xi\in G_\rho$ and $v\in\mcH_\rho$. Define
\begin{equation}\label{Pi}
\Pi(\xi)v=\pi_\rho(\xi)v
\end{equation}
and extend it to the whole $\mcH$ by continuity .
More details about unitarity of $\Pi$ can be found in \cite{29}.

The proof of the following Theorem is inspired by that given in 
\cite{18}.

\begin{theorem} Let $\Pi$ and $\mu$ as before. Then
  \begin{itemize}
    \item
      $\Pi$ is irreducible.
      \item
  Different measures $\mu_1$ and $\mu_2$ on $X$ give rise to inequivalent
  representations $\Pi^1$ and $\Pi^2$ of $\mcG^0$.
  \item
  $\Pi$ is equivalent to $\Xi$.
\end{itemize}
\end{theorem}
\begin{proof}
  A detailed proof of the first two statements can be found in \cite{25}: it can be applied,
  mutatis mutandis, to the case of $PGL(2,\Bbb Q_q)$.

  The third statement requires a computation of matrix coefficents for both
  $\Pi$ and $\Xi$. Fix a partition of $X=\cup_{i=1}^nA_i$ and let
   $z_i=\mu(A_i)$. Let $\xi(x)$ be a function locally
  constant on $A_i$, say $\xi(x)=g_i$ for $x\in A_i$. The matrix coefficent of $\Xi$ with respect to the vacuum vector $\Exp 0$ is
  \begin{align*}
    &\langle\Xi(\xi)\Exp 0,\Exp 0\rangle_{\mbox{\tiny{exp}}} = \\
&q^{{-\frac12\sum_{i=1}^n \mu(A_i)\norm{b(g_i)}^2}}=
q^{{-\frac12\int_X\norm{b(\xi(x))}^2d\mu(x)}}\period
  \end{align*}
  Consider the positive definite function
\begin{equation}
\psi_{z}(g)=q^{-\displaystyle{\frac{z}2\norm{b(g)}^2}}=
q^{\displaystyle{{-z}\dge}+
\displaystyle{ {2zq(1-q^{-\dge})}/{(q^2-1)}}}
\end{equation}
Obvioulsy
$$
\psi_z(g)\geq q^{\displaystyle{{-z}\dge}}\period
  $$

Write the spherical function $\varphi_z$ with respect to the parameter $z$:
\begin{equation}\label{sphfct}
\varphi_z(g)= c(z)q^{-z\dge}+c(1-z)q^{(z-1)\dge}\qquad
\end{equation}

Since
  $0<z_i<\frac12$
  the 
   values $c(z_i)$ are positive while  $c(1-{z_i})$ are negative,
  so that \begin{equation}\label{psiz}
    \psi_{z_i}(g) \geq q^{-z_i{\dge}}=
\frac1{c({z_i})}\varphi_{z_i}(g)
-\frac{c(1-z_i)}{c({z_i})})
q^{(1-z_i)-\dge}\geq\frac1{c({z_i})}\varphi_{z_i}(g)
\end{equation}
   
Hence
\begin{equation}\label{ineqmc}
\begin{aligned}
&\frac{\varphi_{z_1}(g_1)\dots
\varphi_{z_n}(g_n)}{c(z_1)\dots c(z_n)}\leq
 \psi_{z_1}(g_1)\dots\psi_{z_n}(g_n)=\\
&q^{{-\frac12\sum_{i=1}^n m(A_i)\norm{b(g_i)}^2}}=
\langle\Xi(\xi)\Exp 0,\Exp 0\rangle_{\mbox{\tiny{exp}}} 
\period
\end{aligned}
\end{equation}

Let 
\begin{equation}
\one_\rho=\otimes_{i=1}^n\frac{\one_{z_i}}{\sqrt{c(z_i)}}\period
\end{equation}

The above inequality  \eqref{ineqmc} becomes
\begin{equation}
\langle\pi_\rho(\xi)\one_\rho,\one_\rho\rangle\leq
\langle\Xi(\xi)\Exp(0),\Exp(0)\rangle_{\mbox{\tiny{exp}}} =\Psi(\xi)
\end{equation}
So that
 the state $\Psi(\xi)$
dominates 
the positive definite function 
$\langle\pi_\rho(\xi)\one_\rho,\one_\rho\rangle$.
Let $L_\rho$ denote the closed linear span, in $\Exp(H^X)$,
of the vectors $\Xi(\xi)\Exp(0)$ with $\xi\in G_\rho$.
 The map $T_\rho:L_\rho\to \mcH_\rho$ defined by
\begin{equation}
T_\rho(\Exp(0))=
(\one_\rho)
\end{equation}
extends by linearity to a 
unitary equivalence between $\mcH_\rho$ and a subrepresentation
of $L_\rho$. Arguing as in \cite{25} one can see that the inclusion is
compatible with the  structure that gives the limit spaces and conclude that 
there
is an inclusion of the direct limit $\lim_{\rightarrow}\mcH_\rho$
into $\Exp(H^X)$. Since both representations are irreducible for
$\mcG^0$ this inclusion is an equivalence.

\end{proof}

\section{ The group $\mcG$ of bounded measurable currents}

Let
\begin{equation}
X_N=\{v\in\mcT\;:d(v,o)\leq N\}
\end{equation}
and 
\begin{equation}
K_N=\{ g\in G\;: g\cdot v= v\qquad\forall v\in X_N\}\period
\end{equation}

A map $F:X\to G$ is said to be bounded if there exists an integer $N$
and a finite number of cosets $g_0K_{N}\dots g_r K_{N}$ such that
$F(X)\subseteq\cup_{i=0}^r g_i K_{N}$.
 Measurability for $F$
 is defined as usual.

For every $g_1, g_2\in G$: define 
$$\Delta(g_1,g_2)=d(g_1\cdot o,
g_2\cdot o)+\sum_{n=1}^\infty\sum_{d(o,v)=n}\frac{d(g_1\cdot v,g_2\cdot
v)}{(q+1)^{n}[1+d(g_1\cdot v,g_2\cdot v)]}
$$
 where $d(v',v)$ is the tree distance between the vertices
$v'$ and $v$. It is clear that
$\Delta$ is a left invariant  metric 
on $G$ generating the topolgy described in Section 2 .

For $\xi_1$, $\xi_2$ in $\mcG^0$ let
\begin{equation}\label{distanza}
\delta(\xi_1,\xi_2)=\sup_{x\in X}\Delta(\xi_1(x),\xi_2(x))
\period
\end{equation}

The group of 
{\bf bounded measurable currents $\mcG$}
 is the completion of $\mcG^0$ with respect to the metric
 defined by \eqref{distanza}.

The following proposition can be proved by standard arguments: 
\begin{proposition}\label{standard}
Let $F$ be a measurable bounded $G$-valued function on $X$ and $\epsilon>0$.
There exists $\xi\in\mcG^0$ such that $\delta(F,\xi)<\epsilon$.
\end{proposition}

It is convenient to introduce the following notation:
\begin{equation}\label{qua}
\log(\psi(g))=-\frac{1}2\norm{b(g)}^2
=-\dge+\frac{2q}{q^2-1}
(1-q^{-\dge})
\end{equation}
so that 
\begin{equation}
\psi_1(g)=\psi(g)=q^{-\frac{1}2\norm{b(g)}^2}
\period
\end{equation}

The following Theorem guarantees that $\Xi$ can be extended to all $\mcG$:

\begin{theorem}
The positive definite function
\begin{equation}
\langle\Xi(\xi)\Exp(0),\Exp(0)\rangle_{\operatorname{\tiny{exp}}}=
q^{\displaystyle{{-\frac{1}2\int_X\norm{b(\xi(x))}^2d\mu(x)}}}
\end{equation}
can be extended to $\mcG$.
\end{theorem}
\begin{proof}
  The proof exibited in \cite{25} works for $PGL(2,\Bbb Q_q)$ as well:
  \begin{itemize}
\item $x\to\norm{b(F(x))}^2$
  is a bounded measurable function so that the integral $\int_X
  \norm{b(F(x))}^2d\mu(x)$ is convergent.
\item
 since \begin{equation*}
    -\dge\leq\log(\psi(g))\leq-\dge+\frac43\leq\dge\qquad
   \end{equation*} 
 one has
 \begin{equation}
|\log(\psi(g))|\leq\dge\leq\Delta(g,e)
\end{equation}
\item If $\xi_n$ is a Cauchy sequence in $\mcG_0$  we may assume that there
  exists a partition $X=\cup_{j=1}^JA_j$ such that $\xi_n(x)=k_j^n$ and
  $\xi_m(x)=k_j^m$ on $A_j$. Hence 
 \begin{equation}\label{end}
\begin{aligned}
&\left|\sum_{j=1}^J\mu(A_j)(\norm{b(g_j^n)}^2-
\norm{b(g_j^m)}^2)\right|\leq
2\sum_{j=1}^J\mu(A_j)\left|(\log(\psi(g_j^n))-
\log(\psi(g_j^m))\right|=\\
&2\sum_{j=1}^J\mu(A_j)\left|\log\left(\frac{\psi(g_j^n)}
{\psi(g_j^m)}\right)\right|
\period
\end{aligned}
 \end{equation}
\item
  Assuming that $d(o,g_j^n\cdot o)=k_j^n< d(o,g_j^m\cdot o)=k_j^m$ one has 
  \begin{equation}
\frac{\psi(g_j^n)}{\psi(g_j^m)}=
q^{\displaystyle{k_j^m -k_j^n}}
q^{\displaystyle{\frac{2q(q^{-k_j^m}-q^{-k_j^n})}{q^2-1}}}
  \end{equation}
  Hence
\begin{equation}\label{fine}
|\log(\frac{\psi(g_j^n)}{\psi(g_j^m)})|\leq
d(o,g_j^m\cdot o)-d(o,g_j^n\cdot o)\leq d(g_j^m\cdot o,g_j^n\cdot o)
\leq \Delta(g_j^m,g_j^n)\leq\delta(\xi_m,\xi_n) \period
\end{equation}

\item Putting together
  \eqref{end} and \eqref{fine} one sees that $\int_X\norm{b(\xi_n(x))}^2d\mu(x)$
  is a Cauchy sequence.
\item The result now follows by Proposition~\eqref{standard} and standard limit arguments
  (see Lemma 3.5 of \cite{18}). 
   \end{itemize}
  
\end{proof}

\end{document}